\documentclass{article}

\usepackage{amsfonts}
\usepackage{amssymb}
\usepackage{enumerate}
\usepackage{graphicx}

\newtheorem{ttt}{Theorem}[section]
\newtheorem{llll}[ttt]{Lemma}
\newtheorem{ccc}[ttt]{Claim}
\newtheorem{eee}[ttt]{Example}
\newtheorem{sss}[ttt]{Statement}
\newtheorem{ddd}[ttt]{Definition}
\newtheorem{qqq}[ttt]{Question}
\newtheorem{ppp}[ttt]{Problem}
\newtheorem{cccc}[ttt]{Corollary}

\newcommand{\bt}{\begin{ttt}}
\newcommand{\bl}{\begin{llll}}
\newcommand{\bc}{\begin{ccc}}
\newcommand{\bex}{\begin{eee}}
\newcommand{\bs}{\begin{sss}}
\newcommand{\bd}{\begin{ddd} \upshape}
\newcommand{\bq}{\begin{qqq}}
\newcommand{\bprob}{\begin{ppp}}
\newcommand{\bcor}{\begin{cccc}}

\newcommand{\bp}{\noindent\textbf{Proof }}
\newcommand{\br}{\noindent\textbf{Remark }}

\newcommand{\et}{\end{ttt}}
\newcommand{\el}{\end{llll}}
\newcommand{\ec}{\end{ccc}}
\newcommand{\eex}{\end{eee}}
\newcommand{\es}{\end{sss}}
\newcommand{\ed}{\end{ddd}}
\newcommand{\eq}{\end{qqq}}
\newcommand{\eprob}{\end{ppp}}
\newcommand{\ecor}{\end{cccc}}

\newcommand{\ep}{\hspace{\stretch{1}}$\square$\medskip}

\newcommand{\lab}[1]{\label{#1}}

\newcommand{\NN}{\mathbb{N}}

\newcommand{\RR}{\mathbb{R}}
\newcommand{\SSS}{\mathbb{S}}

\newcommand{\ga}{\gamma} 

\newcommand{\e}{\varepsilon}

\newcommand{\bdy}{\partial}

\title{Level Sets of Differentiable Functions of Two Variables with
Non-vanishing Gradient}

\author{M. Elekes}

\begin{document}

\maketitle 

\begin{abstract}
We show that if the gradient of $f:\RR^2\rightarrow\RR$ exists everywhere
and is nowhere zero,
then in a neighbourhood of each of its points the level set
$\{x\in\RR^2:f(x)=c\}$ is homeomorphic either to an open interval or to the
union of finitely many open segments passing through a  point. The second
case holds only at the points of a discrete set. We also investigate the
global structure of the level sets.
\end{abstract}

\insert\footins{\footnotesize{MSC codes: Primary 26B10; Secondary 26B05}}
\insert\footins{\footnotesize{Key Words: Implicit Function Theorem,
Non-vanishing gradient, locally homeomorphic}}

\section*{Introduction}
The Inverse Function Theorem is usually proved under the assumption that
the mapping is continuously differentiable. In \cite{RR} S. Radulescu and
M. Radulescu generalized this theorem to mappings that are only
differentiable, namely they proved that if $f:D\rightarrow\RR^n$ is
differentiable on an open set $D\subset\RR^n$ and the derivative $f'(x)$
is non-singular for every $x\in D$, then $f$ is a local diffeomorphism.

It is therefore natural to ask whether the Implicit Function Theorem, which
is usually derived from the Inverse Function Theorem, can also be proved
under these more general assumptions. (In addition, this question is also
related to the Gradient Problem of C. Weil, see \cite{Qu}, and is motivated
by \cite{EKP} as well, where such a function of two variables with
non-vanishing gradient is
used as a tool to solve a problem of K. Ciesielski.) In \cite{Bu}
Z. Buczolich gave a negative answer to this
question by constructing a differentiable function $f:\RR^2\rightarrow\RR$
of non-vanishing gradient such that $\{x\in\RR^2:f(x)=0\} = \{(x,y):y=x^2
\textrm{ or } y=0 \textrm{ or } y=-x^2\}$. Indeed, this example shows that
the level set is not homeomorphic to an open interval in any neighbourhood
of the origin.

The goal of our paper is to show that such a level set cannot be `much
worse' than that. The main result is a kind of Implicit Function Theorem
(Theorem~\ref{main}), stating that if
$f:\RR^2\rightarrow\RR$ is a differentiable function of non-vanishing
gradient, then in a neighbourhood of each of its points a level set is
homeomorphic either to an open interval or to the union of finitely many
open segments passing through a point. We also show, that the set of points
at which the second case holds has no point of accumulation. 

In order to prove this local result we have to investigate the global
structure of the level set as well. In Section~\ref{secJordan} we apply
the theory of plane continua to prove that the level sets (extended by
$\infty$) consist of Jordan curves (Corollary~\ref{Jordan}). From
this result we show that the level set consists of arcs which have tangents
at each of
their points, only finitely many arcs can
meet at a point and the set of points where these arcs meet has no point
of accumulation. Finally, from all these the above `Implicit Function
Theorem' will follow.

In addition, we show that the notion `arc with tangents' cannot be replaced by
the more natural notion of `differentiable curve with non-zero derivative'.

\section*{Preliminaries}
The usual compactification of the plane by a single point called $\infty$
is denoted by $\SSS^2=\RR^2\cup\{\infty\}$. Throughout the paper
topological notions such as closure, boundary or component, unless
particularly stated, always refer to $\RR^2$. The notations $cl A$, $int A$
and $\partial A$ stand for the closure, interior and boundary of a set $A$,
respectively. The angle of the two vectors in the plane is denoted by
$ang(x,y)$.  The abbreviations $\{f=c\}$, $\{f<c\}$ etc. stand for
$\{x\in\RR^2 : f(x)=c\}$, $\{x\in\RR^2 : f(x)<c\}$, etc.,
respectively. $B(x,\e)$ is the open disc $\{y\in\RR^2 : |y-x|<\e\}$. The
circle of center $x$ and radius $\e$, is denoted by $S(x,\e)$. By an
\textsl{arc} or a \textsl{Jordan curve} we mean a continuous and injective
function to the plane (or to $\SSS^2$) defined on a closed interval or on a
circle, respectively. (We often do not distinguish between the image of the
function and the function itself.) The \textsl{contingent} of a set
$H\subset\RR^2$ at a point $x\in H$ is the union of those half-lines $L$
that can be written as $L=lim \ L_n$, where $L_n$ is a half-line starting
from $x$ and passing through $x_n\in H$, and $x_n$ is converging to $x$
($x_n\neq x$). (By $L=lim \ L_n$ we mean that the direction of the
half-lines converges.)  We say that $H\subset\RR^2$ has a \textsl{tangent}
(half-tangent) at $x\in H$ if the contingent of $H$ at $x$ is a line
(half-line). A \textsl{continuum} is a compact connected set.

\section{The results}

\bd\lab{nice}
Let $f:\RR^2\rightarrow\RR$ be a differentiable function of non-vanishing
gradient and let $c\in\RR$ be arbitrary. Then $x\in\RR^2$ is called a
\textsl{branching point} of $\{f=c\}$ if it is in the closure of at least
three different components of $\{f\neq c\}$.  We call $H\subset\RR^2$ a
\textsl{nice curve} if $H=\ga\setminus\{\infty\}$ where
\begin{itemize}
\item[(i)] $\ga$ is either an arc in $\RR^2$ between two branching points
or an arc in $\SSS^2$ between a branching point and $\infty$ or a Jordan
curve in $\SSS^2$ containing $\infty$,
\item[(ii)] if $x\in\ga$ and $x$ is not an endpoint, then $\ga$ has a
tangent at $x$, and if $x\in\ga\setminus\{\infty\}$ and $x$ is an endpoint,
then $\ga$ has a half-tangent at $x$.
\end{itemize}
\ed 

\bigskip

\bt\lab{main} Let $f:\RR^2\rightarrow\RR$ be a differentiable function of
non-vanishing gradient and let $c\in\RR$ be arbitrary. Then the set of
branching points has no point of accumulation, and the level set $\{f=c\}$
is the disjoint union (except from the endpoints) of nice curves. Moreover,
if $x\in\RR^2$ is not a branching point, then there exists a neighbourhood
$U$ of $x$ such that $\{f=c\}\cap U$ is homeomorphic to an open interval,
while if $x\in\RR^2$ is a branching point, then there exists a
neighbourhood $U$ of $x$ such that $\{f=c\}\cap U$ is homeomorphic to the
union of finitely many open segments passing through a point.
\et

In the rest of the paper we prove this theorem. As we have already
mentioned in the Introduction, first we examine certain global
properties of the level sets in the next section, and then apply these
results in the last section to obtain Theorem~\ref{main}.

Throughout the proof we assume that $f:\RR^2\rightarrow\RR$ is 
a differentiable function of non-vanishing gradient, $D$ is a component of
$\{f\neq c\}$ and $C$ is a component of $\bdy D$. We can clearly assume
that $c=0$.

\section{The level set consists of Jordan curves}\lab{secJordan}

We start with a lemma that we shall frequently use in the sequel.

\bl\lab{contingent} Let $a<b$ and $c<d$ be real numbers and $F\subset
[a,b]\times [c,d]$ be a closed set that has infinitely many points on each
vertical line that meets the rectangle. Then there is a point in $F$ at
which the contingent of $F$ is not contained in a line.
\el

\bp
Suppose, on the contrary, that at every point of $F$ the contingent of $F$
is contained in a line.
Our assumption is that for every $a\leq x_0\leq b$ the set $\{(x,y)\in F :
x=x_0\}$ is infinite. Therefore we can choose a point of accumulation of
this set for every $a\leq x_0\leq b$, and thus we obtain a function
$g:[a,b]\to [c,d]$. As $F$ is closed, $(x_0,g(x_0))\in F$ and because of
the way the point was chosen, the
contingent of $F$ at this point must be contained in the vertical
line. Hence for every $a\leq x_0\leq b$
\begin{equation}\label{difference}
{lim}_{x\to x_0} \left|\frac{g(x)-g(x_0)}{x-x_0}\right| = +\infty.
\end{equation}
But this is impossible, as \cite[IX. 4. 4]{Sa} states that for any function
of a real variable the set of points at which (\ref{difference}) holds is
of measure zero.
\ep

The proof of the next lemma is a straightforward calculation, so
we omit it.

\bl\lab{tangent} For every $x\in\{f=0\}$ the contingent of $\{f=0\}$ at $x$
is contained in the line perpendicular to $f'(x)$.
\el

\begin{figure}[!ht]
\begin{center}
\includegraphics[scale=.8]{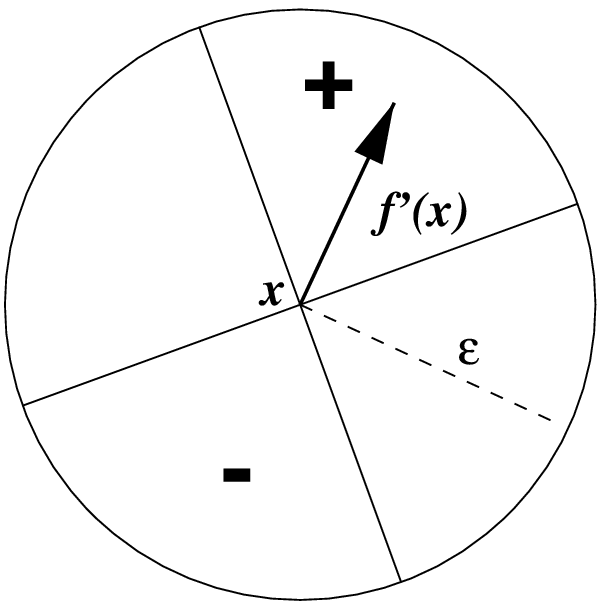}
\caption{}\label{figcross}
\end{center}
\end{figure}

An easy consequence is the following (see Figure~\ref{figcross}).

\bl\lab{cross} For every $x\in\{f=0\}$ there exists an $\e>0$ such that if
$y\neq x$ and $y\in B(x,\e)$, then $\ ang(f'(x),y-x)\leq\frac{\pi}{4}$
implies that $f(y)>0$, while $ang(\textrm{-}f'(x),y-x)\leq\frac{\pi}{4}$
implies that $f(y)<0$.
\el

The next lemma is the basic tool in the proof of the fact that $C \cup
\{\infty\}$ is a Jordan curve (Corollary~\ref{Jordan}, the main result of
the present section).

\bl\lab{lc}\
\begin{enumerate}[(i)]
\item $\{f=0\}$ is locally connected.
\item $\bdy D$ has no bounded component.
\item $\bdy D\cup \{\infty\}$ is connected.
\item Both $\bdy D$ and $\bdy D\cup \{\infty\}$ are locally connected.
\item Both $C$ and $C\cup \{\infty\}$ are connected and locally connected.
\end{enumerate}
\el

\bp
\begin{enumerate}[(i)]
\item \lab{one} If $\{f=0\}$ is not locally connected at $x\in\{f=0\}$,
then for some $\e>0$ we can find a sequence $K_n\ (n\in\NN)$ of distinct
components of $\{f=0\}\cap cl{B(x,\e)}$ converging (with respect to the
Hausdorff metric) to a continuum $K$ such that $x\in K$ and $K_n\cap
K=\emptyset$ for every $n\in\NN$ (\cite[I. 12. 1]{Wh} asserts that if the
compact set $M\subset\RR^2$ is not locally connected at a point $m$, then
for some $\e>0$ there exists a sequence $M_n$ of components of $M\cap
cl{B(x,\e)}$ converging to a continuum $N$ such that $m\in N$ and $M_n\cap
N=\emptyset$ for every $n\in\NN$). We claim that $K_n\cap S(x,\e)\neq
\emptyset$ for every $n\in\NN$. Indeed, if this is the case, then there
exists a Jordan curve inside ${B(x,\e)}$ that encloses $K_n$ and that is
disjoint from $\{f=0\}$ (\cite[VI. 3. 11]{Wh} states that if $M$ is a
component of a compact set $N\subset\RR^2$, then there exists a Jordan
curve in the $\e$-neighbourhood of $M$ that encloses $M$ and is disjoint
from $N$). But then the sign of
$f$ is constant on the curve and thus $f$ attains a local extremum inside
the curve, which contradicts the assumption that the gradient nowhere
vanishes.

Let us now divide $S(x,\e)$ into three sub-arcs of equal length. At least
one of these pieces must intersect infinitely many of the sets $K_n\
(n\in\NN)$. Let us call this subsequence $K_{n_i}\ (i\in\NN)$. This chosen
sub-arc
can be separated from $x$ by a narrow rectangle (see Figure~\ref{separate}). 

\begin{figure}[!ht]
\begin{center}
\includegraphics[scale=.8]{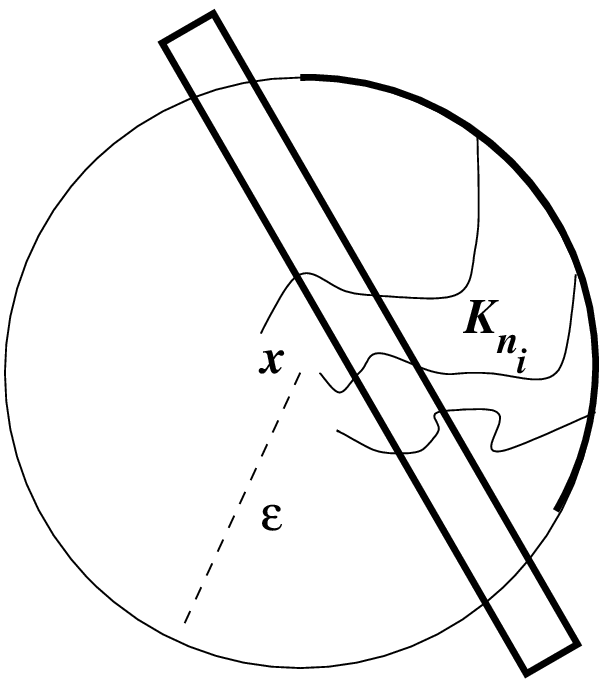}
\caption{}\label{separate}
\end{center}
\end{figure}

As
$K_{n_i}\rightarrow K$, $x\in K$, $K_{n_i}\cap S(x,\e)\neq \emptyset$ and
$K_{n_i}$ is connected $(i\in\NN)$, for all $i\in\NN$ large enough,
$K_{n_i}$ must `cross' the rectangle from the sub-arc to a point close to
$x$. As the sets $K_{n_i}\ (i\in\NN)$ are disjoint, we obtain infinitely
many different points of intersection on the lines considered in
Lemma~\ref{contingent}, so we can
apply this lemma (to a rotated copy of our rectangle) with
$F=\{f=0\}$. But then we get that the contingent of $\{f=0\}$ is not
contained in a line, which contradicts Lemma~\ref{tangent}.

\item\lab{two} Suppose that a component $C$ of $\bdy D$ is contained in
$B((0,0),R)$. Applying \cite[VI. 3. 11]{Wh} (see (\ref{one})) to $\bdy
D\cap cl{B((0,0),R)}$
we can again find a Jordan curve $\ga$ around the component $C$ that is
disjoint from $\bdy D$. So $\ga$ is either contained in $D$ or disjoint
from $D$. In the first case the sign of $f$ is constant on the curve, so
$f$ attains a local extremum, which is impossible. In the second case $D$
must be inside $\ga$, as $D$ is connected and at least one point of it is
inside $\ga$, since there are points of $\bdy D$ inside $\ga$. But no
component of $\{f=0\}$ can be bounded, since it would result in a local
extremum
of $f$ as $f$ vanishes on the boundary of the components.

\item\lab{three} Easily follows from (\ref{two}).

\item\lab{four} If $x\in\RR^2$, then we can repeat the argument of
(\ref{one}). The only difference is that we prove $K_n\cap S(x,\e)\neq
\emptyset\ (n\in\NN)$ by applying (\ref{two}).

In the case of $\infty$ we argue as follows. Let $\pi : \SSS^2 \setminus
\{s\} \to \RR^2$ be the usual stereo-graphic projection from a point $s \in
\SSS^2$ which is not in $\bdy D\cup \{\infty\}$. As this is a local
homeomorphism in a neighbourhood of $\infty$, it is sufficient to prove
that $\pi(\bdy D\cup \{\infty\})$ is locally connected at the point
$x=\pi(\infty)\in\RR^2$. Suppose, on the contrary, that this is not
true. It is enough to find a point $y\neq x$ at which the contingent
of $\pi(\bdy D\cup \{\infty\})$ is not contained in a line, since $\pi$ is
a local diffeomorphism at $p=\pi^{\textrm{-}1}(y)$, and so in this case the
contingent of $\bdy D\cup \{\infty\}$ at $p$ is also not contained in a
line, which contradicts Lemma~\ref{tangent}. But we can again obtain this
by the argument of (\ref{one}), once we show that $K_n\cap
S(x,\e)\neq\emptyset\ (n\in\NN)$. So let $K_n \subset B(x,\e)$. As $x
\notin K_n$, its inverse image $\pi^{\textrm{-}1}(K_n)$ is bounded, and
moreover $K_n\cap
S(x,\e)\neq\emptyset$, therefore $\pi^{\textrm{-}1}(K_n)$ is
a component of $\bdy D\cup \{\infty\}$. (Note that $K_n$ was
originally a component of $\pi(\bdy D\cup \{\infty\}) \cap
cl{B(x,\e)}$.) But this is impossible by (\ref{two}).

\item $C$ is clearly connected and by (\ref{two}) unbounded. Thus $C \cup
\{ \infty \}$ is connected as well. What remains to prove is that these
sets are locally connected. If $x \in \RR^2$, then we can repeat the
argument of (\ref{one}). The only difference is that $K_n \cap S(x,\e) \neq
\emptyset$ simply follows from the connectedness of $C$. For the case of
$\infty$
we repeat the proof of (\ref{four}). We again simply use the connectedness
of $C$ to verify that $K_n \cap S(x,\e) \neq \emptyset$.
\end{enumerate}
\ep

\bcor\lab{accessible} Every $x \in \bdy D \cup \{ \infty \}$ can be
accessed from any point of $D$ by an arc in $D \cup \{ x \}$.
\ecor

\bp \cite[\S 61, II. 11]{Ku} states that if $E$ is a component of the
complement of a locally connected and closed set $F$, then every point of
$\bdy E$ can be accessed from $E$ by a connected set. A remark at the end
of the proof adds that it can also be accessed by an arc.

Suppose first that $x\neq\infty$.
Let us apply this theorem to $E=D$ and $F=\bdy D$. Indeed, $D$ is a
component of of the complement of $\bdy D$, which is closed and by
(\ref{four}) of Lemma~\ref{lc} locally connected. Thus the statement follows.

To see that the case $x=\infty$ is similar, note first that $\bdy D \cup \{
\infty \}$ is also locally connected by (\ref{four}) of
Lemma~\ref{lc}. Choose a point $s\in \bdy D$ and let $\pi : \SSS^2 \setminus
\{s\} \to \RR^2$ be the usual stereo-graphic projection from the point
$s$. We can now apply the above result to the open set $\pi(D)\subset\RR^2$
and its boundary $\bdy\pi(D) = \pi(\bdy D \cup \{\infty\} \setminus
\{s\})$, since this latter set is clearly locally connected. Therefore we
obtain an arc $\ga\subset \pi(D)$ accessing
$\pi(\infty)$, and then $\pi^{\textrm{-}1}(\ga)$ is the required arc.
\ep

\bs\lab{twoarcs} For every $x \in \bdy D$ there exists a Jordan curve in
$\bdy D \cup \{ \infty \}$ containing both $x$ and $\infty$.
\es

Before the proof of this statement we need a technical lemma. It is surely
well known, but we could not find it in the literature, so we
include a proof here.

\bl Let $\ga_1$ and $\ga_2$ be arcs between a point $x \in \RR^2$ and $\{
\infty \}$ such that they are disjoint except from the endpoints, and let
$G$ be one of the two components of the complement of the Jordan curve
formed by the two arcs. In addition, let $\varphi$ be a Jordan curve in
$\RR^2$ which contains $x$ in its interior. Then there exists a sub-arc of
$\varphi$ that joins $\ga_1$ and $\ga_2$ such that $\varphi \subset G$
except from the endpoints.
\el

\begin{figure}[!ht]
\begin{center}
\includegraphics[scale=.8]{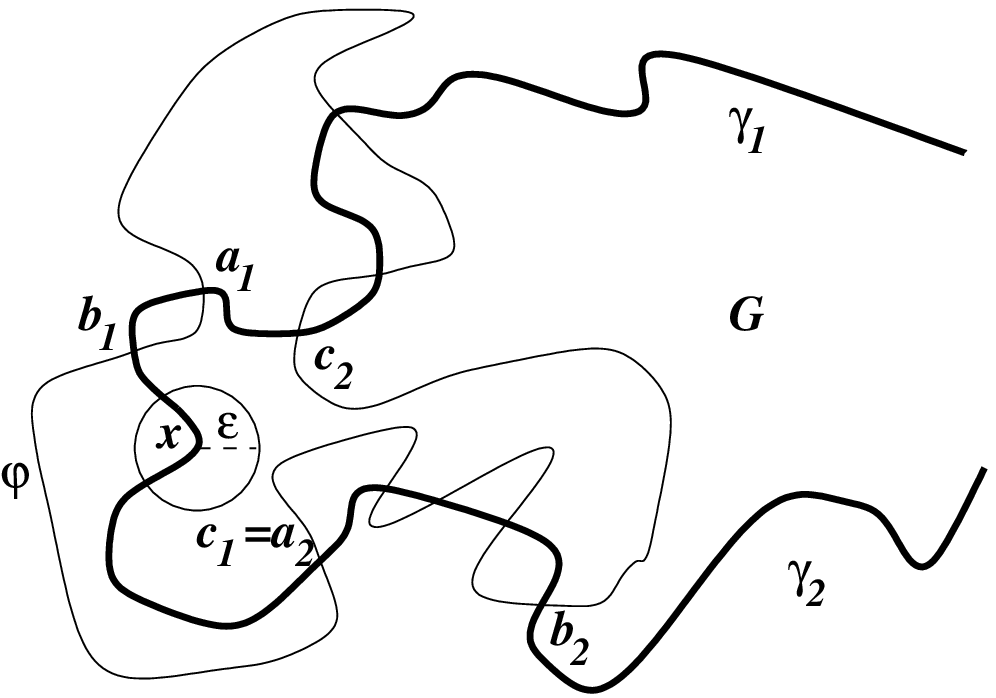}
\caption{}\label{sub-arc}
\end{center}
\end{figure}

\bp $\varphi : [0,1] \rightarrow \RR^2 \setminus B(x,\e)$ for some $\e>0$,
and we may assume that $\varphi(0)=\varphi(1)$ is on one of the arcs. Let
$c_0=0$ and let us construct inductively a sequence of triples in
the following way (see Figure~\ref{sub-arc}). 

For every $n=1,2,\dots$ put
$a_n=c_{n-1}$ and define $c_n$ as the smallest number (greater than
$a_n$) for which $\varphi(c_n)$ is already on the other arc than
$\varphi(a_n)$. Finally, let $b_n$ be the largest number (smaller
than $c_n$) for which $\varphi(b_n)$ is still on the same arc as
$\varphi(a_n)$. We claim that for some $n\in\NN$ $\ c_n=1$. Otherwise, $a_n\
(n\in\NN)$ is a strictly increasing sequence converging to a number $a\in
[0,1]$. So the sequence $[b_n,c_n]\ (n\in\NN)$ of disjoint intervals also
converges to $a$, thus $\varphi([b_n,c_n])$ converges to the point
$\varphi(a)$. But these images are arcs between $\ga_1$ and $\ga_2$, and
they are disjoint from $B(x,\e)$, hence they cannot converge to a point. 

If we replace the sub-arcs of $\varphi$ between $\varphi(a_i)$ and $\varphi(b_i) \
(i=1, \dots, n)$ by the corresponding sub-arcs of $\ga_1$ or $\ga_2$, we
obtain a continuous closed curve $\psi$ of the same rotation number around
$x$ as $\varphi$, but the rotation number around $x$ of $\varphi$ is 1 (or -1),
thus $\psi$ must intersect $G$. But then for the corresponding $[b_i,c_i]$,
$\varphi|_{[b_i,c_i]}$ is the required sub-arc.
\ep

Now we prove Statement~\ref{twoarcs}.

\medskip

\bp
By Corollary~\ref{accessible} there exists an arc $\ga_1$ from $x$ to
$\infty$ in $D\cup\{x,\infty\}$. As $f$ has no local extremum, there
must be another component $E$ of $\{f\neq 0\}$ (of different sign than $D$)
such that $x \in \bdy E$. Let $\ga_2$ be a similar arc in
$E\cup\{x,\infty\}$ from $x$ to $\infty$.

Since $\bdy D \cup \{\infty\} \subset \SSS^2$ is closed and locally
connected, it is locally arc-wise connected (this is \cite[\S 50,
II. 1]{Ku}). Moreover, it is connected by (\ref{three}) of Lemma~\ref{lc},
hence arc-wise connected. Therefore
there exists an arc $\psi_1$ in $\bdy D \cup \{\infty\}$ from $x$ to
$\infty$. $\ga_1$ and $\ga_2$ joined together is a Jordan curve in
$\SSS^2$, thus it splits its complement into two components, and $\psi_1$
must be contained in one of them (except from the endpoints). Let us denote
the other component by $G$. It is sufficient to show that there exists an arc
$\psi_2$ between $x$ and $\infty$ in $(G \cap \bdy D) \cup
\{x,\infty\}$. First we check that this set is locally connected. We only
have to consider $x$ and $\infty$. But at these points the locally
connected set $\bdy D \cup \{\infty\}$ is cut into two pieces by a
Jordan curve that intersects the set in no additional point (in a
neighbourhood), hence the set must be locally connected on `both sides' of
the Jordan curve. So  $(G \cap \bdy D) \cup \{x,\infty\}$ is locally
connected. Now we show that $x$ and $\infty$ are in the same component of
$(G \cap \bdy D) \cup \{x,\infty\}$. Otherwise, the component of $x$
must be bounded, thus it can be enclosed by a Jordan curve $\varphi$ that is
disjoint from $(G \cap \bdy D) \cup \{x,\infty\}$. But then by the
previous lemma we obtain an arc from $\ga_1$ (which is in $D$) to $\ga_2$
(which is not in $D$) such that this arc is in $G$ except from its
endpoints. This is a contradiction. 

The common component of $x$ and $\infty$ in $(G \cap \bdy D) \cup
\{x,\infty\}$ is locally connected, as it is a component of a locally
connected set. By the same argument as at the beginning of the previous
paragraph we obtain that it is arc-wise
connected. Hence we can construct an arc $\psi_2$ in $(G \cap \bdy D) \cup
\{x,\infty\}$ between $x$ and $\infty$, and then $\psi_1$ and $\psi_2$
together form the required Jordan curve.
\ep

Now we are able prove the main result of this section.

\bcor\lab{Jordan}
Let $f:\RR^2\rightarrow\RR$ be a differentiable function of non-vanishing
gradient, $D$ be a component of $\{f\neq 0\}$ and $C$ be a component of
$\bdy D$. Then $C \cup \{\infty\}$ is a Jordan curve.
\ecor

\bp
We apply \cite[\S 52, VI. 1]{Ku}, which asserts that if a locally connected
continuum (of at least two points) contains no $\theta$-curve and no
separating point, then it is a Jordan curve. (A $\theta$-curve is the union
of three arcs between two points such that the arcs are disjoint except
from the endpoints. A point is a separating point if its complement is
disconnected.) 

By (\ref{two}) of Lemma~\ref{lc} we know
that $C \cup \{\infty\}$ is a locally connected continuum. First we check
that it contains no $\theta$-curve. The complement (in $\SSS^2$) of a
$\theta$-curve consists of three Jordan domains, and $D$ must be contained
in one of them. But the boundary of this domain is the union of two arcs of
the $\theta$-curve, therefore the third arc cannot be in $\bdy D$,
which is a contradiction. We still have to check that $C \cup \{\infty\}$
contains no separating point, which easily follows from
Statement~\ref{twoarcs}.
\ep

\section{The rest of the proof of Theorem~\ref{main}}

\bs\lab{finitely}
For every $x\in\RR^2$ and $R>0$ there are only finitely many components
of $\{f\neq 0\}$ intersecting the disc $B(x,R)$.
\es

\begin{figure}[!ht]
\begin{center}
\includegraphics[scale=.8]{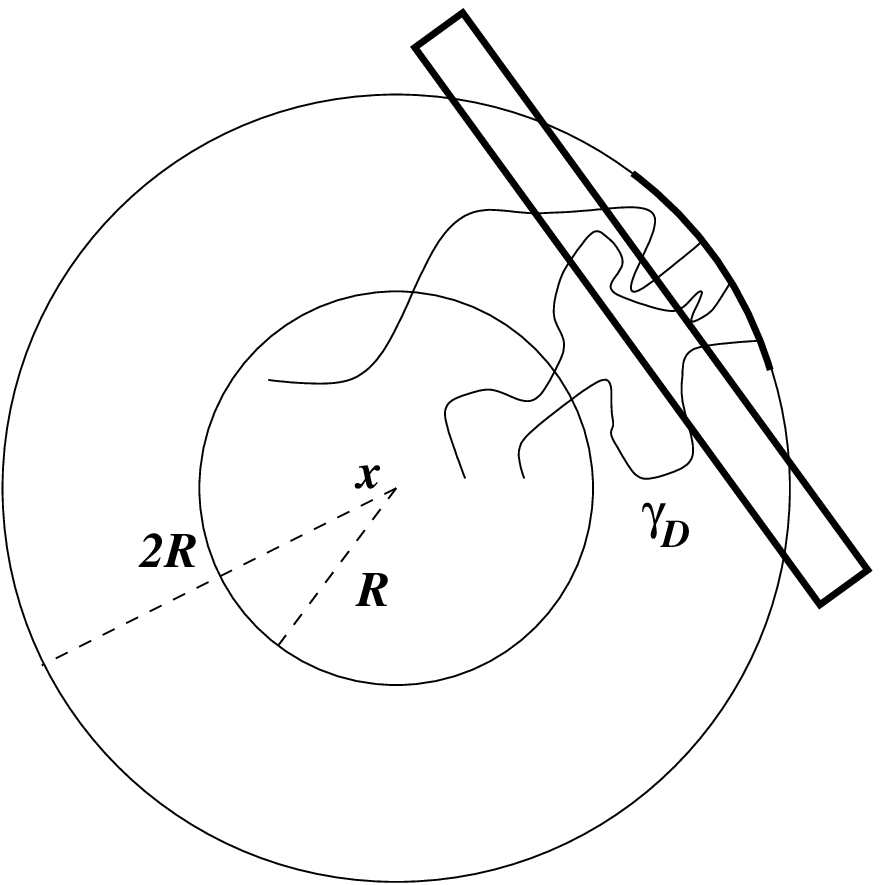}
\caption{}\label{separate2}
\end{center}
\end{figure}

\bp
Suppose that the converse is true.
As every component is unbounded (otherwise we could find a local extremum),
there exists an
arc $\ga_D \subset D$ for every component $D$ such that $\ga_D$ joins
$S(x,R)$ and $S(x,2R)$. We may assume that $\ga_D \subset cl B(x,2R)$. Indeed,
denote by $x_D$
the first point along the arc $\ga_D$ that is on $S(x,2R)$ and cut off the
rest of $\ga_D$. Our aim is to
get a contradiction by a similar argument as  in the proof of
Lemma~\ref{lc}. So choose a small sub-arc of $S(x,2R)$ (e.g. one tenth of
the circle) that contains infinitely many of the points $x_D$, and
separate this sub-arc from $S(x,R)$ by a narrow rectangle as in
Figure~\ref{separate2}. 

Let us consider a segment in our rectangle as in Lemma~\ref{contingent}. It
is intersected by all of the disjoint arcs $\ga_D$. Hence we obtain
infinitely many points on our segment that are in distinct components of
$\{f\neq 0\}$, therefore there must also be infinitely many points of the set
$\{f=0\}$ on the segment. Thus Lemma~\ref{contingent} can be applied, and
we get a contradiction by Lemma~\ref{tangent}.
\ep

\bcor\lab{accumulation}
The set of branching points has no point of accumulation in $\RR^2$.
\ecor

\bp
By the previous statement there are only finitely many components of
$\{f\neq 0\}$ intersecting $B(x,1)$. Therefore the proof is complete once
we show that for two distinct components $D_1$ and $D_2$ there are at most two
branching points in $\bdy D_1 \cap \bdy D_2$. Suppose, on the contrary,
that $x_1, x_2, x_3\in\RR^2$ are three such points, and let $d_1\in D_1,\ 
d_2\in D_2$ be arbitrary. By Lemma~\ref{accessible} we can join $x_1, x_2$
and $x_3$ by six 
disjoint arcs (except from the endpoints) to $d_1$ and $d_2$, and so we obtain
three disjoint arcs (except from the endpoints) from $d_1$ to $d_2$ through the
points $x_1, x_2$ and $x_3$. Let us denote them by $\ga_1, \ga_2$ and
$\ga_3$. One of these arc, say $\ga_2$, must be surrounded by the Jordan
curve formed by the other two arcs. But then $x_2$ cannot be a branching
point, since if it were on the boundary of a third component of
$\{f\neq 0\}$, then this component would be inside the Jordan curve, hence
bounded, which is impossible.
\ep

We need one more lemma to prove the main result, our version of the Implicit
Function Theorem (Statement~\ref{homeomorphic}).

\bl\lab{intersection}
Let $D_1$ and $D_2$ be two components of $\{f\neq 0\}$. Then $\bdy D_1 \cap
\bdy D_2 = C_1 \cap C_2$, where  $C_1$ and $C_2$ are
components of $\bdy D_1$ and $\bdy D_2$, respectively. Moreover, this
intersection is a sub-arc (possibly empty) of both of the curves $C_1$ and
$C_2$. An endpoint of such a sub-arc is either a branching point or
$\infty$, but the other points of the sub-arc are not branching points.
\el

\bp
In order to justify the first equality it is sufficient to show that $\bdy
D_1$ cannot intersect two components of $\bdy D_2$. Suppose, on the contrary,
that $C_2$ and $C_2'$ are two such components. By Corollary~\ref{Jordan}
they are Jordan curves in $\SSS^2$ (apart from $\infty$), hence they split
$\RR^2$ into three domains $G_1,G_2$ and $G_3$ such that the boundaries of
these domains are $C_2, C_2'$ and $C_2 \cup C_2'$,
respectively. Consequently, $D_1 \subset G_3$ and $D_2 \subset G_3$ must
hold. As $\bdy D_1$ intersects both $C_2$ and $C_2'$, by Lemma
\ref{accessible} there exists an arc $\ga$ in $D_2$ (apart from the
endpoints) which joins $C_2$ and $C_2'$. But $\ga$ splits $G_3$ into two
parts and $D_1$ must be contained in one of them, which contradicts
e.g. $C_2 \subset \bdy D_1$.

To prove that the intersection is a sub-arc of e.g. $C_1$, let us denote by
$\ga_{D_1}$ and $\ga_{D_2}$
two Jordan curves in $\SSS^2$ such that $\ga_{D_1} \setminus \{\infty\} =
C_1$ and $\ga_{D_2} \setminus \{\infty\} = C_2$. It is sufficient to show
that there are no three points $x,y$ and $z$ on $\ga_{D_1} \setminus
\{\infty\}$ in this order such that $x,z\in \ga_{D_2}$ and $y \notin
\ga_{D_2}$. Suppose, on the contrary, that there are three such
points. On the sub-arc of $\ga_{D_1} \setminus \{\infty\}$ from $y$ to $x$
there exists a first point $x'$ in $\ga_{D_1} \cap \ga_{D_2}$. Similarly,
the first point of this intersection on $\ga_{D_1} \setminus \{\infty\}$ in
the other direction is denoted by $z'$. But these two points, $x'$ and $z'$
are
connected by a sub-arc of $\ga_{D_1} \setminus \{\infty\}$ and a sub-arc of
$\ga_{D_2} \setminus \{\infty\}$, which form a Jordan curve in $\{f=0\}$,
which results in a local extremum of $f$, a contradiction. 

Now we have to check that if an endpoint of a sub-arc is not $\infty$, then
it is a branching point. If $x\in\RR^2$ is such an endpoint, then $x$ and
$\infty$ splits $\ga_{D_1}$ and $\ga_{D_2}$ into sub-arcs between $x$ and
$\infty$ such that three of these arcs are disjoint except from their
endpoints. But these three arcs divide the plane into three domains,
therefore at
least one of them is disjoint from both $D_1$ and $D_2$. As $x$ is on the
boundary of all three domains, it must meet the closure of a component
different from $D_1$ and $D_2$. Thus $x$ is a branching point.

We now show that (except for the endpoints) no point of the sub-arc can be a
branching point. Suppose, on the contrary, that there exists such a point
$y$. By Lemma~\ref{accessible} we can join two other points of the sub-arc
to two points $d_1\in D_1$ and $d_2\in D_2$ such that the obtained Jordan
curve contains $y$ in its interior. As $y$ is a branching point, a
third component $D_3$ must intersect the interior of this Jordan curve,
moreover, it cannot intersect the curve itself, therefore it must be
enclosed by the Jordan curve. But then $D_3$ is bounded, and $f$ attains a
local extremum inside it, which is impossible.
\ep

Now we can prove the most important statement of Theorem~\ref{main}.

\bs\lab{homeomorphic}
If $x\in\RR^2$ is not a branching point, then there exists a neighbourhood
$U$ of $x$ such that $\{f=0\}\cap U$ is homeomorphic to an open interval,
while if $x\in\RR^2$ is a branching point, then there exists a
neighbourhood $U$ of $x$ such that $\{f=0\}\cap U$ is homeomorphic to the
union of finitely many open segments passing through a point.
\es

\bp
Let $x\in\RR^2$ be arbitrary. By Statement~\ref{finitely} there are only
finitely many components of $\{f\neq0\}$ intersecting $B(x,1)$, therefore
for some $\e>0$ the point $x$ is on the boundary of every component that
intersects $B(x,\e)$. We can also assume by Corollary~\ref{accumulation}
that $B(x,\e)$ contains no branching point, with the possible exception of
$x$ itself. 

Let us now first suppose that $x$ is not a branching point, that is only
two components of $\{f\neq0\}$ intersect $B(x,\e)$. As $f$ has no
local extremum, every point of $\{f=0\}$ is on the boundary of both
components. Thus by Lemma~\ref{intersection} we obtain that
$\{f=0\} \cap B(x,\e)$ is the intersection of an arc with $B(x,\e)$ such
that the endpoints of the arc are outside the disc. If we now choose an
open neighbourhood $U$ of $x$ inside $B(x,\e)$ by (\ref{one}) of Lemma
\ref{lc} such that $U \cap  \{f=0\}$ is connected, then this intersection
must be homeomorphic to an open interval.

Let us now consider the case when $x$ is a branching point. As $f$ has no
local extrema, every point in $\{f=0\} \cap  B(x,\e)$ is on the common
boundary of at least two components of $f \neq 0$. Since $x$ is the
only branching point in the disc, we obtain that $\{f=0\} \cap  B(x,\e)$ is
the intersection of $B(x,\e)$ and the disjoint union (except from the point
$x$) of finitely many arcs starting from $x$ and running to branching
points outside the disc (indeed, we apply Lemma~\ref{intersection} to every
pair of components intersecting the disc). If we now choose an open
neighbourhood $U$ of $x$ inside $B(x,\e)$ such that $U \cap  \{f=0\}$ is
connected, then this intersection must be homeomorphic to a half-open
interval on each of the above arcs starting from $x$. So the only thing
that remains to show is that there is an even number of these arcs, which
easily follows from the fact that the components of $\{f \neq 0\}$ surrounding
$x$ must be of alternating signs.
\ep

Our next goal is to complete the proof of Theorem~\ref{main}.

\bd
Let $D$ be a component of $\{f \neq 0\}$, $C$ be a component of $\bdy D$ and
$\ga$ be a Jordan curve such that $\ga = C\cup\{\infty\}$. Then $\ga$ is
separated into sub-arcs by the branching points and $\infty$. We call these
sub-arcs (together with their endpoints) \textsl{edges}.
\ed

\bs
The edges are disjoint except from their endpoints. If $\varphi$ is an edge,
then $\varphi \setminus \{\infty\}$ is a nice curve (see Definition
\ref{nice}).
\es

\bp
If two edges correspond to the same component of $\{f \neq 0\}$, then they are
clearly disjoint except from their endpoints. Let now $\varphi_1 \subset \ga_1
\subset \bdy D_1$ and  $\varphi_2 \subset \ga_2 \subset \bdy D_2$ be two edges
corresponding to the distinct components $D_1$ and $D_2$ such that they
have a point $x$ in common which is not an endpoint of at least one of
them. As an edge is an arc between branching points and $\infty$, no point
inside an edge can be a branching point, thus $x$ is not a branching
point. By Lemma~\ref{intersection} $\ (\ga_1 \setminus \{\infty\}) \cap
(\ga_2 \setminus \{\infty\})$ is a sub-arc of $\ga_1$ such that its
endpoints are branching points or $\infty$ and the other points of the
sub-arc are not branching points, therefore this sub-arc must agree with
$\varphi_1 \setminus \{\infty\}$. Similarly $(\ga_1 \setminus \{\infty\}) \cap
(\ga_2 \setminus \{\infty\})$ must agree with $\varphi_2 \setminus
\{\infty\}$, therefore the two edges coincide.

To show that $\varphi \setminus \{\infty\}$ is a nice curve note that it
clearly satisfies (i) of Definition~\ref{nice} (by Corollary
\ref{Jordan}). To see that (ii) is also satisfied, let $x\in\varphi \setminus
\{\infty\}$. By Lemma~\ref{tangent} the contingent of $\varphi$ at $x$ is at
most the line perpendicular to $f'(x)$. 

Suppose first that $x$ is an endpoint. An easy compactness argument shows
that the contingent of $\varphi$ at $x$ cannot be empty, thus it is sufficient
to show that it does not contain both possible half-lines. If we apply
Lemma~\ref{cross} to $x$ we obtain two opposite sectors of $B(x,\e)$
containing the directions of the two possible half-lines, and moreover
$\varphi \cap B(x,\e)$
must also be contained in these two sectors. But $\varphi$ is an arc which
starts from $x$ and never returns there, from which easily follows that it
cannot approach $x$ arbitrarily close inside both sectors. 

Suppose next that $x$ is not an endpoint of $\varphi$, thus it is not a
branching point. We have to show that the contingent contains both possible
half-lines. By Statement~\ref{homeomorphic} we can find a neighbourhood $U$
of $x$ such that $\{f=0\} \cap U$ is a sub-arc of $\varphi$. Then Lemma
\ref{cross} provides a small disc $B(x,\e)$ inside $U$ which consists of four
sectors, two opposite sectors $S_1$ and $S_3$ containing the directions of
the two possible half-lines and two other sectors $S_2$ and $S_4$, in which
$f$ is positive and negative, respectively. Because of the different signs,
these latter sectors cannot be connected by an arc in $\{f\neq 0\}$,
therefore for every $\delta<\e$ there must be points of $\{f=0\}$ on both
semi-circles of center $x$ and of radius $\delta$ running from $S_2$ to
$S_4$ and crossing
through $S_1$ or $S_3$. Consequently, the contingent contains both
half-lines by an easy compactness argument.
\ep

Therefore the proof of Theorem~\ref{main} is complete.

\medskip

\begin{figure}[!ht]
\begin{center}
\includegraphics[scale=.8]{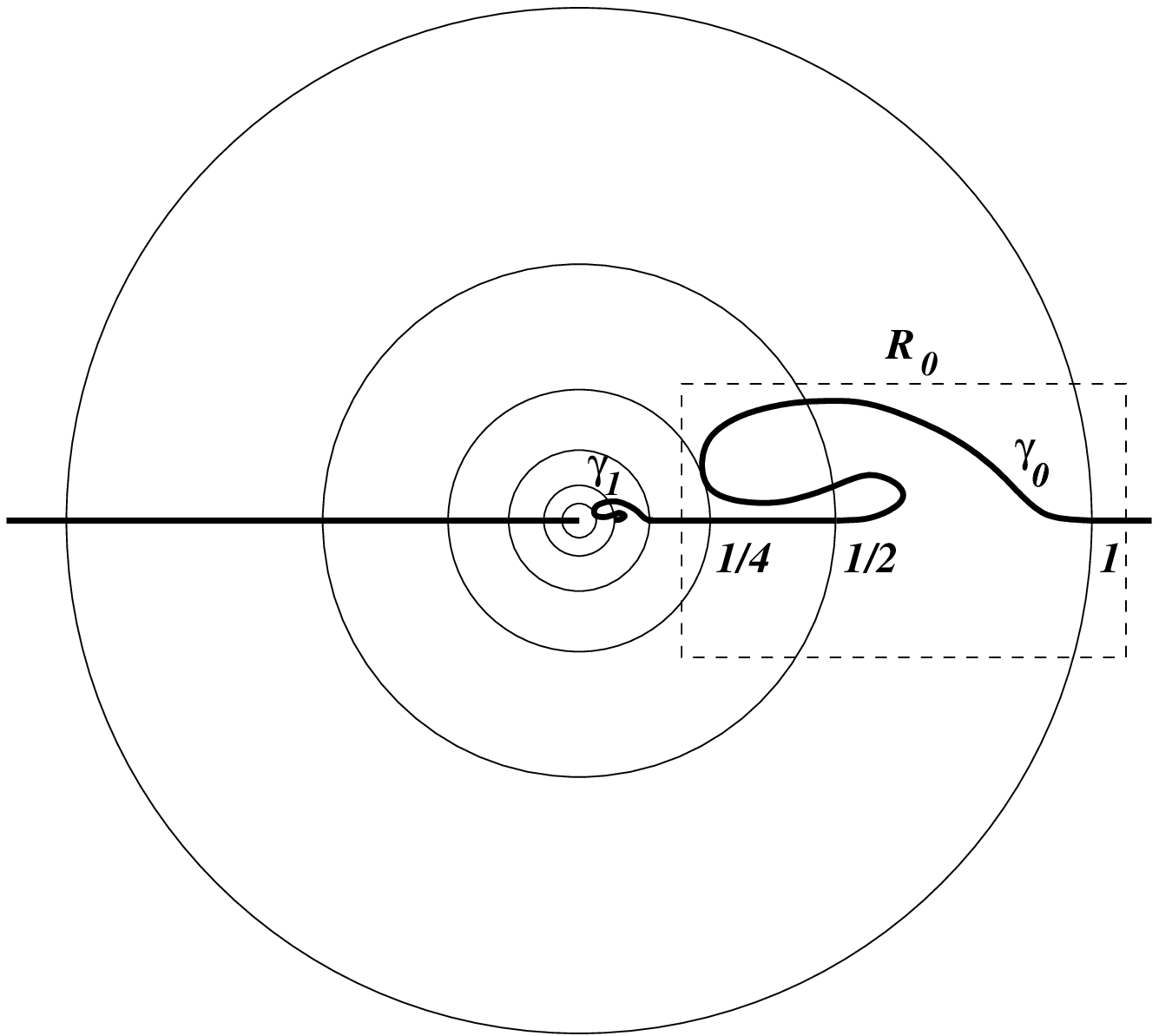}
\caption{}\label{gamma}
\end{center}
\end{figure}

\br
As we have already mentioned in the Introduction, it is not true that the level
sets consist of differentiable arcs with non-zero derivative: 

We can argue
as follows. For every $n\in\NN$ let $\ga_n$ be a smooth curve in the
(closed) region bounded by $y=0,\ y=x^2,\ S(0,\frac{1}{2^{3n}})$ and
$S(0,\frac{1}{2^{3n+2}})$ such that its endpoints are
$(0,\frac{1}{2^{3n+1}})$ and $(0,\frac{1}{2^{3n}})$, it meets
$S(0,\frac{1}{2^{3n+2}})$ and such that if we continue $\ga_n$ by two
horizontal segments to the left and to right, then we get a smooth curve
(see Figure~\ref{gamma}). 
Let $H\subset\RR^2$ be the union of the curves
$\ga_n\ (n\in\NN)$, the horizontal segments connecting them and the two
half-lines $\{(x,y):\ x\leq 0,\ y=0\}$ and $\{(x,y):\ x\geq 1,\ y=0\}$. One
can show that the Jordan curve $H \cup \{\infty\}$ is the level set
$\{f=0\}$ of a suitable differentiable
function $f:\RR^2\rightarrow \RR$ of non-vanishing gradient, but it is not
a differentiable curve of non-zero derivative. We only sketch the proofs
here. 

To construct the function $f$ let us fix for every $n\in\NN$ a rectangle
$R_n$ (as in Figure~\ref{gamma}) containing $\ga_n$ in its interior, such
that $R_n$ is between the parabolas $y=2x^2$ and $y=\textrm{-}2x^2$. It is
not hard to see that there exists a diffeomorphism $\Phi^{(n)} : R_n \to
R_n$ for
which $\Phi^{(n)}(x,y) \in \{y=0\}$ iff $(x,y) \in \ga\cap R_n$ (this
latter set is the
union of $\ga_n$ and the two horizontal segments). We can also clearly
assume that $\Phi^{(n)}$ coincides with the identity function close to the edges
of $R_n$. Now define
\[
f(x,y) = 
\left\{
\begin{array}{ll}
\Phi^{(n)}_2 (x,y) & \textrm{if } (x,y)\in R_n,\\
y                 & \textrm{otherwise},\\
\end{array}
\right.
\]
where $\Phi^{(n)}_2$ is the second coordinate function of
$\Phi^{(n)}=(\Phi^{(n)}_1,\Phi^{(n)}_2)$. It is easy to check that 
$f'(0,0)$ exists and is equal to $(0,1)$, and one can also see that the
function $f$ satisfies all the other requirements.

To see that the Jordan curve
$H\cup\{\infty\}$ cannot have a non-zero derivative at the origin, we
parametrize it by a function $\varphi$ such that $\varphi(0)=(0,0)$, and
consider
\[
\left|\frac{\varphi(t)-\varphi(0)}{t-0}\right| =
\left|\frac{\varphi(t)}{t}\right|.
\]
Note that there exists a sequence $t_n \rightarrow 0 \ (t_n>0)$ such that
$\varphi(t_n) = (0,\frac{1}{2^{3n+1}})$ and an other sequence $t_n'
\rightarrow 0 \ (t_n'>0)$ for which we have $\varphi(t_n') \in
S(0,\frac{1}{2^{3n+2}})$, while $t_n' \geq t_n$. This shows that the right
hand side derivative of $\varphi$ at $0$ cannot be a finite, non-zero value.

\medskip

Finally we pose the following problem.

\bprob
Characterize the level sets $\{f=c\}$ of differentiable functions $f:\RR^2 \to
\RR$ of non-vanishing gradient. 
\eprob

\end{document}